\documentclass[%
preprint,
bibnotes,
amsmath,
amssymb,
aomart
showkeys
]{revtex4-1}

\usepackage{graphicx}
\usepackage{dcolumn}
\usepackage{bm}
\usepackage{amsmath}

\begin{document}

\def\mean#1{\left< #1 \right>}
\title{A linear transformation to accelerate the convergence of the negative binomial series}
\author{Liborio I. Costa}
\affiliation{8405 Winterthur, Switzerland}%
\date{February 05, 2017}

\begin{abstract}
A linear sequence transformation is defined that accelerates the convergence of the negative binomial series when the terms of the binomial have the same sign.
The transformed series can be used to extend the region of applicability of the Taylor expansion of ln$(1+x)$ and to compute the incomplete beta function.
\end{abstract}


\maketitle


                                                \section{Introduction}  
The binomial theorem is one of the fundamental theorems in mathematics. It states that, for $r, x, y \in \mathbb{R}$, and $|x/y|<1$ \cite{Graham}:

\begin{equation}\label{Binomial}
\left(x+y\right)^r =  \sum\limits_{k=0}^{\infty} \frac{(r)_k}{k!}x^ky^{r-k}
\end{equation}

The identity in Eq.~\ref{Binomial} is widely used in several fields of physics and engineering to expand powers of binomials into more manageable polynomials.
On the other hand, when $r$ is negative and $x$ and $y$ have the same sign, the series on the right-hand side of Eq.~\ref{Binomial} converges
rather slowly when the values of $x$ and $y$ are of comparable magnitude. As $x/y \rightarrow 1$, the convergence becomes so slow that accumulation of round-off errors arising from the summation of terms with higher indexes can severely affect the final result. Moreover, it does not converge at all when $x/y=1$ and $r\leq-1$. Under such conditions, in fact, the binomial theorem is notoriously not valid and the identity in Eq.~\ref{Binomial} does not hold. As a consequence, the use of the polynomial expansion in Eq.~\ref{Binomial} for negative exponents becomes practically inconvenient when $x=O(y)$, and can not be used when $x=y$. The convergence properties of slowly convergent series can be improved by applying series transformation techniques \cite{Weniger}. Although nonlinear transformations
are in general expected to be more effective than linear transformations, the latter have the advantage of being easier to apply and can, in some cases, be very effective too \cite{Michalski,Mosig,Zaiger}.
Accordingly, in this work a linear sequence transformation is introduced to improve the convergence properties of the binomial series when $r<0$ and $x$ and $y$ have the same sign. Numerical examples indicate that the rate of convergence of the transformed series is with good approximation not affected by the ratio $x/y$, resulting in a significant speed-up when $x=O(y)$. Notably, the modified series converges to the value of a binomial with negative exponent also when $x/y=1$, i.e. when the binomial theorem is not valid.

                                                \section{Definitions}  
\indent \it{Definition} \rm1. Given the power of a binomial $\left(x+y\right)^r$, let us define the associated $0^{th}$-level partial sum as:

\begin{equation}\label{sum_0_level}
s^0_{xyr}(n) =  \sum\limits_{k=0}^{n} a_k  \qquad n \in\mathbb{N}_0
\end{equation}

where $a_k$ is the element $k$ of the infinite binomial series in Eq.~\ref{Binomial}:

\begin{equation}\label{a_k}
a_k =  \frac{(r)_k}{k!}x^ky^{r-k}
\end{equation}

and $(r)_k$ is the falling factorial.\\

\indent \it{Definition} \rm2. Given the power of a binomial $\left(x+y\right)^r$, let us define the associated $j^{th}$-level partial sum as:

\begin{equation}\label{sum_j_level}
s^j_{xyr}(n) =  \frac{1}{2}s^{j-1}_{xyr}(n) + \frac{1}{4}s^{j-1}_{xyr}(n-1) + \frac{1}{4}s^{j-1}_{xyr}(n+1) \qquad n \geq j; n,j \in\mathbb{N}
\end{equation}\\

Thus, with Eq.~\ref{sum_0_level} and Eq.~\ref{sum_j_level} we have introduced the elements of the sequences of partial sums $S^0_{xyr}(n)=\{s^0_{xyr}(n)\}$ and $S^j_{xyr}(n)=\{s^j_{xyr}(n)\}$, respectively. Clearly, $S^0_{xyr}(n)$ is the sequence of partial sums of the infinite binomial series.\\

\indent \it{Remark} \rm1. From the two definitions above, it follows that the sequence of the first elements of the sequences $S^j_{xyr}(n)$, $S^j_{xyr}(j)=\{s^j_{xyr}(j)\}$, are calculated using the first $2j+1$ elements of the sequence $S^0_{xyr}(n)$, as illustrated in Tab. \ref{Table1}, i.e. to calculate $s^j_{xyr}(j)$, only the first $2j+1$ terms $a_k$ are required.\\

\begin{table}
\caption{Example of the dependency between $s^j(j)$ and $s^0(n)$ for $j=3$. The elements marked in bold are those required to calculate $s^3(3)$, cf. Definition 2.}
\centering
\begin{tabular} { l | c c c  c c c }
\hline
     &  $j$  &  0  &    1    &  2   &  3  \\
 $n$   &     &     &         &      &     \\
\hline
 0   &     & \boldmath $s^0(0)$  &    &        &         \\
 1   &     & \boldmath $s^0(1)$  &   \boldmath $s^1(1)$ &        &        \\
 2   &     & \boldmath $s^0(2)$  &   \boldmath $s^1(2)$ & \boldmath $s^2(2)$       &         \\
 3   &     & \boldmath $s^0(3)$  &   \boldmath $s^1(3)$ & \boldmath $s^2(3)$       &  \boldmath $s^3(3)$  \\
 4   &     & \boldmath $s^0(4)$  &   \boldmath $s^1(4)$ & \boldmath $s^2(4)$       &   $s^3(4)$  \\
 5   &     & \boldmath $s^0(5)$  &   \boldmath $s^1(5)$ &     $s^2(5)$             &   $s^3(5)$  \\
 6   &     & \boldmath $s^0(6)$  &             $s^1(6)$ &     $s^2(6)$             &   $s^3(6)$  \\
\hline
\end{tabular}
\label{Table1}
\end{table}

\indent \it{Remark} \rm2. By definition, the sequence $S^j_{xyr}(j)$ is a linear transformation of the sequence
$S^0_{xyr}(n)$, and therefore one can write:
\begin{equation}\label{modified_expansion}
s^j_{xyr}(j) =  \sum\limits_{k=0}^{2j} c_{kj}a_k
\end{equation}
where the coefficients $c_{kj} \geq 0$ are uniquely determined for any $j$, and can be calculated once for all recursively by applying Definition 1 and Definition 2.\\

As an example, let us consider the element $s^1_{xyr}(1)$:
\begin{equation}\label{example_S11}
\begin{aligned}
s^1_{xyr}(1) {} & =  \frac{1}{2}s^{0}_{xyr}(1) + \frac{1}{4}s^{0}_{xyr}(0) + \frac{1}{4}s^{0}_{xyr}(2) = \\
                &  = \frac{1}{2}(a_0+a_1) + \frac{1}{4}a_0 + \frac{1}{4}(a_0+a_1+a_2) = a_0 + \frac{3}{4}a_1 + \frac{1}{4}a_2
\end{aligned}
\end{equation}

By comparing Eq. \ref{example_S11} and Eq. \ref{modified_expansion} it is easy to see that:  $c_{0,1}=1, c_{1,1}=3/4, c_{2,1}=1/4$.\\

In the following section it will be shown that the sequence of partial sums $S^j(j)$, which was introduced by the simple linear transformation defined by Eq.~\ref{sum_j_level},
converges significantly faster than the slowly convergent negative binomial series as $x/y \rightarrow 1$.
In section IV also some applications of the method will be illustrated.

                                                \section{Numerical examples}  
In this section numerical results are reported to show the convergence properties of the sequence of partial sums $S^j(j)$. It is worth keeping in mind
that we limit our analysis to the case $r<0$ and $0\leq x/y \leq 1$. Under such conditions, defined $ q=x/y $, we have that:

\begin{equation}\label{Binomial_with_q}
(x+y)^r=y^r(q+1)^r
\end{equation}

from which it follows that:
\begin{equation}\label{Sequence_with_q}
s^j_{xyr}(j)=y^rs^j_{q1r}(j)
\end{equation}

Therefore, for the purpose of comparing $S^j(j)$ with the original sequence of partial sums $S^0(n)$, it suffices to study the expansion of $(q+1)^r$ with $0\leq q \leq 1$, which is done in the following.

The columns of Tables \ref{Table2} to \ref{Table5} report the elements of $S^0_{q1r}(n)$ and $S^j_{q1r}(n)$ for the case $r=-1$, while the elements of $S^j_{q1r}(j)$ are on the diagonals. Each table corresponds to a different value of $q$, increasing from $q=0.1$ in Table \ref{Table2}, to $q=1$ in Table \ref{Table5}. In each table, the first element of $S^0_{q1r}(n)$ and of $S^j_{q1r}(j)$ which is equal to $(q+1)^r$ to the first 6 decimal digits is marked in bold.
When $q$ is at least an order of magnitude smaller than 1, as in Table \ref{Table2}, $S^0_{q1r}(n)$ converges faster than $S^j_{q1r}(j)$, but both do it after a few terms. Thus, using one or the other sequence does not change much.
The picture changes drastically when $q=O(1)$, as in Tables \ref{Table3} to \ref{Table5}. As $q$ approaches unity, more and more terms are required to let $S^0_{q1r}(n)$ converge, and not even 90 terms
are enough when $q=0.9$ (cf. Table \ref{Table4}). Differently, $S^j_{q1r}(j)$ converges after a few terms in all cases. The extreme case is reached for $q=1$ in Table \ref{Table5}: $S^j_{q1r}(j)$ converges immediately for $j=1$, in contrast to the non-convergent alternating series $S^0_{q1r}(n)$.
Figures 1-3 show the same type of data for $r$ negative non-integer ($r=-0.5$ in Figure 1, $r=-3\sqrt[]{2}$ in Figure 2) and negative integer ($r=-10$ in Figure 3) in terms of the relative truncation error:

\begin{equation}\label{Error_S0}
\text{Relative Truncation Error} = \frac{|s^0_{q1r}(n)-(q+1)^r|}{(q+1)^r}
\end{equation}
and
\begin{equation}\label{Error_Sj}
\text{Relative Truncation Error} = \frac{|s^j_{q1r}(j)-(q+1)^r|}{(q+1)^r}
\end{equation}

respectively. The results plotted in the figures confirm the same trend observed for the case $r=-1$: $S^j_{q1r}(j)$ converges within a few terms irrespective of the values of $r$ and $q$: the relative truncation error falls below $10^{-6}$ for $j=O(\text{ceil}(|r|)$, even when $q=1$ and $r<-1$, despite under such conditions the binomial series is wildly divergent (cf. Figure 2 and Figure 3). This represents a quite surprising result when considering the simplicity of the transformation (cf. Eq. \ref{modified_expansion}) and the fact that neither complex algorithms nor any estimate of the residuals have been necessary, as instead it is the case when using other types of transformations \cite{Weniger,Michalski,Mosig}.

\begin{table}
\caption{Convergence of the sequences for $r=-1$ and $q=0.1$. The correct value is $(q+1)^r=0.909091$. All values are rounded to 6 decimal digits. The subscript $q1r$ is implicit.}
\centering
\begin{tabular} { l | c c c c c c c c c }
\hline
 $n$ &  $S^0(n)$    &  $S^1(n)$     &  $S^2(n)$     &  $S^3(n)$  &  $S^4(n)$ &  $S^5(n)$ &  $S^6(n)$ &  $S^7(n)$ &  $S^8(n)$  \\
\hline
  0  &	1.000000	&	    -	    &	-        	&	-        &	-        &	   -     &	    -    &	   -     &	-  \\
  1  &	0.900000	&	0.927500	&	-        	&	-        &	-        &	   -     &	    -    &	   -     &	-  \\
  2  &	0.910000	&	0.907250	&	0.912819	&	-        &	-        &	   -     &	    -    &	   -     &	-  \\
  3  &	0.909000	&	0.909275	&	0.908718	&	0.909846 &	-        &	   -     &	    -    &	   -     &	-  \\
  4  &	0.909100	&	0.909073	&	0.909128	&	0.909015 &	0.909244 &	   -     &	    -    &	   -     &	-  \\
  5  &	0.909090	&	0.909093	&	0.909087	&	0.909098 &	0.909076 &	0.909122 &	    -    &	   -     &	-  \\
  6  &	\textbf{0.909091}	&	0.909091	&	0.909091	&	0.909090 &	0.909092 &	0.909088     &	0.909097 &	- &	-  \\
  7  &	0.909091	&	0.909091	&	0.909091	&	0.909091 &	0.909091 &	0.909091 &	0.909091 &	0.909092 &	-  \\
  8  &	0.909091	&	0.909091	&	0.909091	&	0.909091 &	0.909091 &	0.909091 &	0.909091 &	0.909091 &	\textbf{0.909091}  \\
  9  &	0.909091	&	0.909091	&	0.909091	&	0.909091 &	0.909091 &	0.909091 &	0.909091 &	0.909091 &	0.909091  \\
  10 &	0.909091	&	0.909091	&	0.909091	&	0.909091 &	0.909091 &	0.909091 &	0.909091 &	0.909091 &	0.909091  \\
\end{tabular}
\label{Table2}
\end{table}

\begin{table}
\caption{Convergence of the sequences for $r=-1$ and $q=0.5$. The correct value is $(q+1)^r=0.666667$. All values are rounded to 6 decimal digits. The subscript $q1r$ is implicit.}
\centering
\begin{tabular} { l | c c c c c c}
\hline
 $n$ &  $S^0(n)$    &  $S^1(n)$     &  $S^2(n)$     &  $S^3(n)$  &  $S^4(n)$ &  $S^5(n)$  \\
\hline
  0 &	1.000000	&	-	&	-	&	- &	- &	-  \\
  1 &	0.500000  	&	0.687500	&	-	&	- &	- &	- \\
  2 &	0.750000	&	0.656250	&	0.667969	&	- &	- &	-  \\
  3 &	0.625000  	&	0.671875	&	0.666016	&	0.666748 &	- &	- \\
  4 &	0.687500	&	0.664063	&	0.666992	&	0.666626 &	0.666672 &	-  \\
  5 &	0.656250	&	0.667969	&	0.666504	&	0.666687 &	0.666664 &	\textbf{0.666667}  \\
  6 &	0.671875	&	0.666016	&	0.666748	&	0.666656 &	0.666668 &	0.666667  \\
  7 &	0.664063	&	0.666992	&	0.666626	&	0.666672 &	0.666666 &	0.666667  \\
  8 &	0.667969	&	0.666504	&	0.666687	&	0.666664 &	0.666667 &	0.666667  \\
  9 &	0.666016	&	0.666748	&	0.666656	&	0.666668 &	0.666667 &	0.666667  \\
  10 &	0.666992	&	0.666626	&	0.666672	&	0.666666 &	0.666667 &	0.666667  \\
  11 &	0.666504	&	0.666687	&	0.666664	&	0.666667 &	0.666667 &	0.666667  \\
  12 &	0.666748	&	0.666656	&	0.666668	&	0.666667 &	0.666667 &	0.666667  \\
  13 &	0.666626	&	0.666672	&	0.666666	&	0.666667 &	0.666667 &	0.666667  \\
  14 &	0.666687	&	0.666664	&	0.666667	&	0.666667 &	0.666667 &	0.666667  \\
  15 &	0.666656	&	0.666668	&	0.666667	&	0.666667 &	0.666667 &	0.666667  \\
  16 &	0.666672	&	0.666666	&	0.666667	&	0.666667 &	0.666667 &	0.666667  \\
  17 &	0.666664	&	0.666667	&	0.666667	&	0.666667 &	0.666667 &	0.666667  \\
  18 &	0.666668	&	0.666667	&	0.666667	&	0.666667 &	0.666667 &	0.666667  \\
  19 &	0.666666	&	0.666667	&	0.666667	&	0.666667 &	0.666667 &	0.666667  \\
  20 &	\textbf{0.666667}	&	0.666667	&	0.666667	&	0.666667 &	0.666667 &	0.666667  \\
\end{tabular}
\label{Table3}
\end{table}

\begin{table}
\caption{Convergence of the sequences for $r=-1$ and $q=0.9$. The correct value is $(q+1)^r=0.526316$. All values are rounded to 6 decimal digits. The subscript $q1r$ is implicit.}
\centering
\begin{tabular} { l | c c c c  }
\hline
 $n$ &  $S^0(n)$    &  $S^1(n)$     &  $S^2(n)$     &  $S^3(n)$   \\
\hline
  0 &	1.000000	&	-	        &	-	        &	-       \\
  1 &	0.100000  	&	0.527500	&	-	        &	-       \\
  2 &	0.910000	&	0.525250	&	0.526319	&	-       \\
  3 &	0.181000  	&	0.527275	&	0.526313	&	\textbf{0.526316}  \\
  4 &	0.837100	&	0.525453	&	0.526318	&	0.526316 \\
  5 &	0.246610	&	0.527093	&	0.526314	&	0.526316 \\
  6 &	0.778051	&	0.525617	&	0.526318	&	0.526316 \\
... &	...      	&	...      	&	...     	&	...      \\
15  &	0.428788    	&	...      	&	...     	&	...      \\
... &	...      	&	...      	&	...     	&	...      \\
30  &	0.546396    	&	...      	&	...     	&	...      \\
... &	...      	&	...      	&	...     	&	...      \\
50  &	0.528757    	&	...      	&	...     	&	...      \\
... &	...      	&	...      	&	...     	&	...      \\
90  &	0.526352    	&	...      	&	...     	&	...      \\
\end{tabular}
\label{Table4}
\end{table}

\begin{table}
\caption{Convergence of the sequences for $r=-1$ and $q=1$. The correct value is $(q+1)^r=0.5$. All values are rounded to 6 decimal digits. The subscript $q1r$ is implicit.}
\centering
\begin{tabular} { l | c c c   }
\hline
 $n$ &  $S^0(n)$    &  $S^1(n)$     &  $S^2(n)$      \\
\hline
  0 &	1.000000	&	-	        &	-	              \\
  1 &	0.000000  	&	\textbf{0.500000}	&	-	             \\
  2 &	1.000000	&	0.500000    &	0.500000           \\
  3 &	0.000000  	&	0.500000	&	0.500000           \\
  4 &	1.000000	&	0.500000    &	0.500000           \\
  5 &	0.000000  	&	0.500000	&	0.500000           \\
\end{tabular}
\label{Table5}
\end{table}

\begin{center}
\begin{figure}
\includegraphics[scale=0.6]{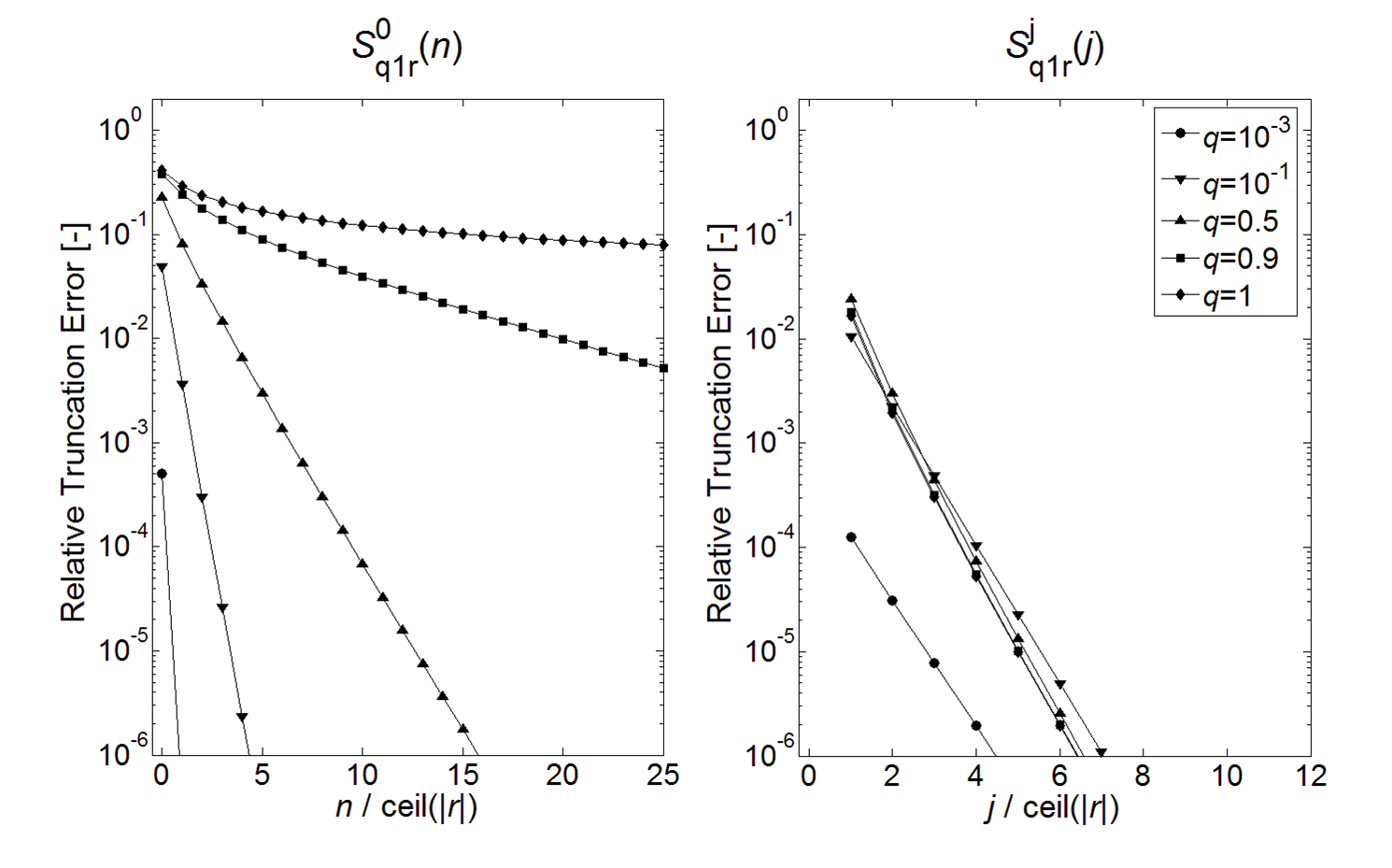}
\caption{\label{Truncation_Error_r_05} Comparison of Relative Truncation Error for $r = -0.5$ and $q$ as in the legend. The symbols of $S^j(j)$ (right figure) for $q=0.9$ and for $q=1$ are not well visible because they are almost superimposed. Lines are a guide to the eye.}
\end{figure}
\end{center}

\begin{center}
\begin{figure}
\includegraphics[scale=0.6]{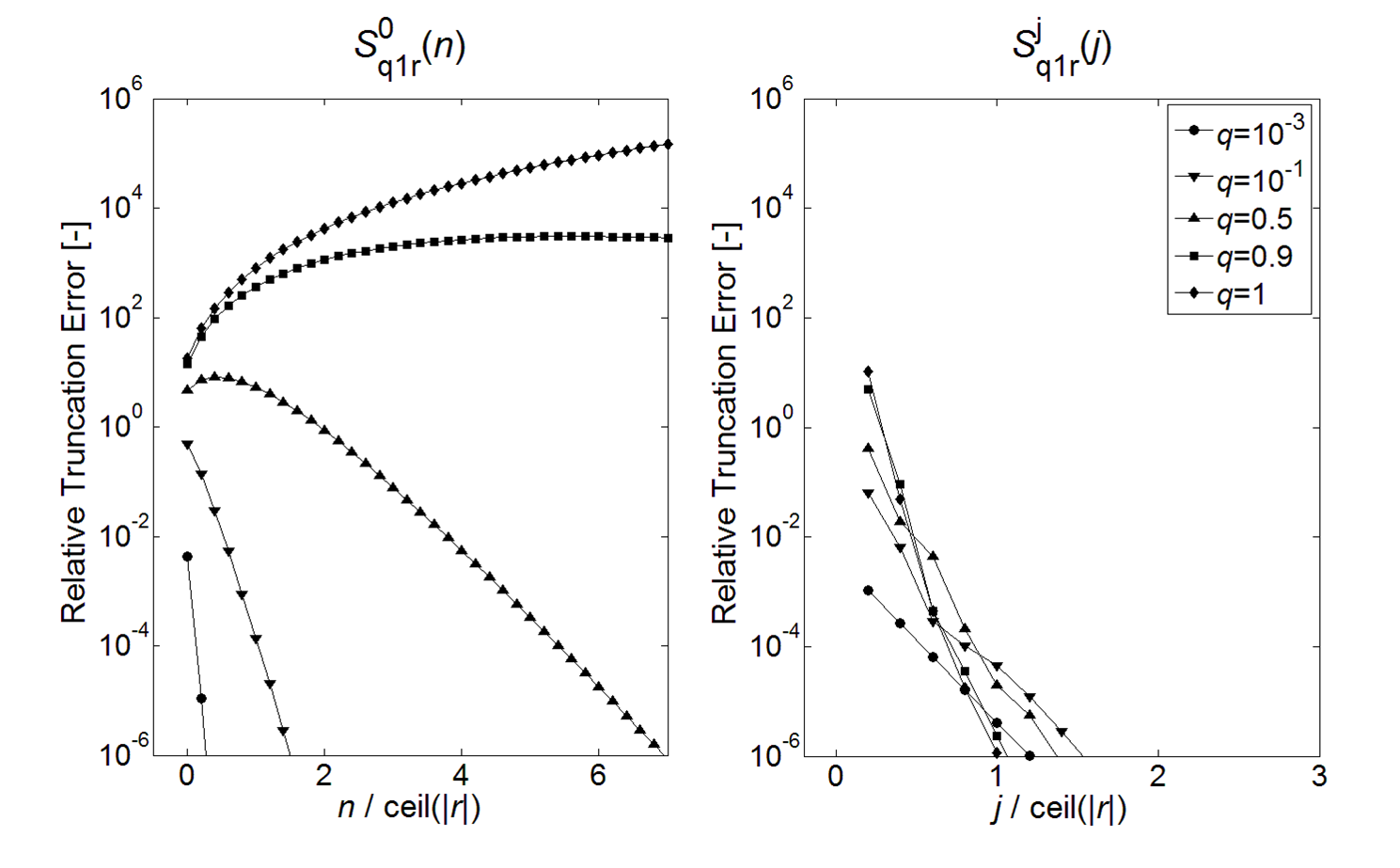}
\caption{\label{Truncation_Error_r_05} Comparison of Relative Truncation Error for $r=-3\sqrt[]{2}$ and $q$ as in the legend. Lines are a guide
to the eye.}
\end{figure}
\end{center}

\begin{center}
\begin{figure}
\includegraphics[scale=0.6]{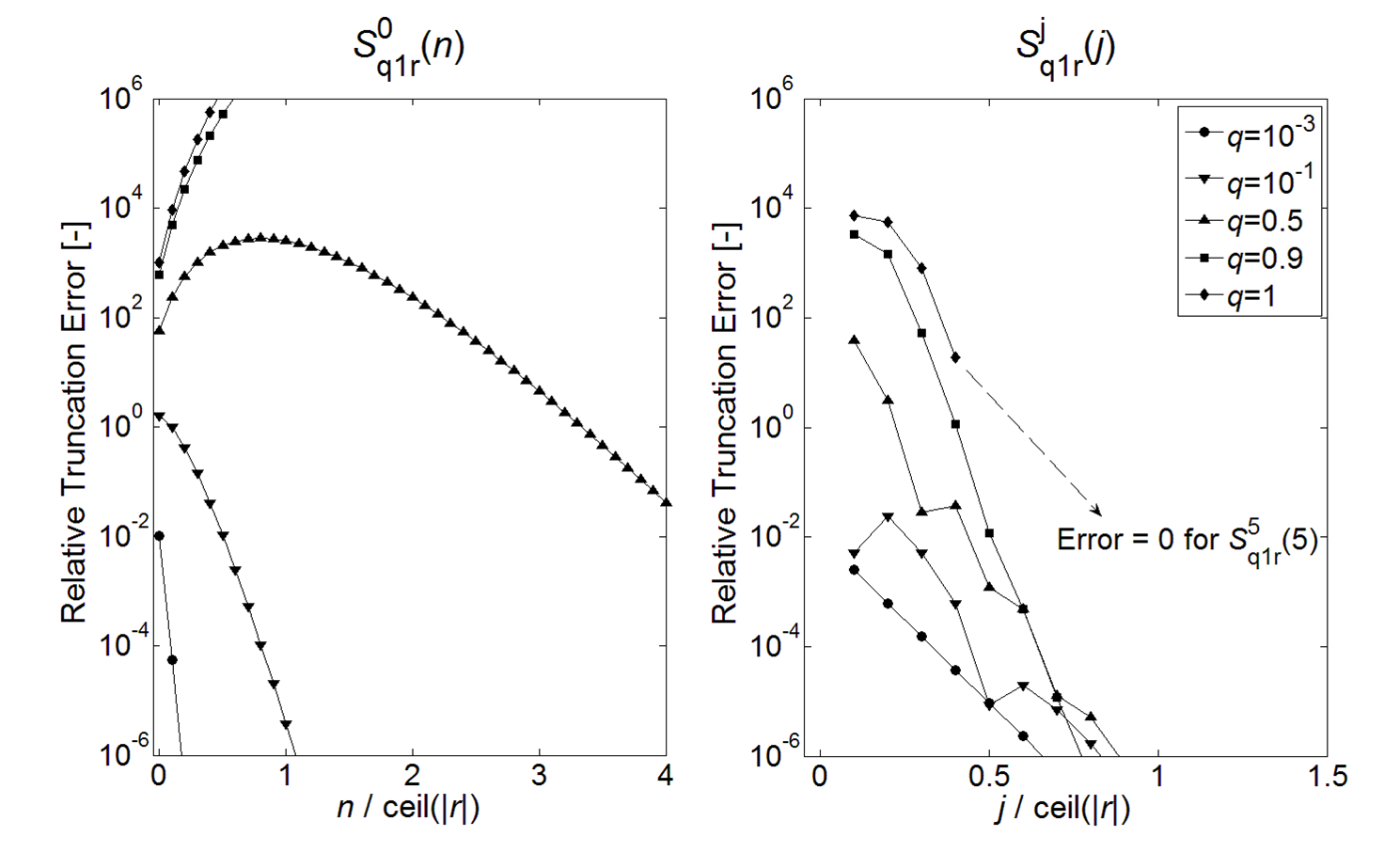}
\caption{\label{Truncation_Error_r_05} Comparison of Relative Truncation Error for $r=-10$ and $q$ as in the legend. Note that only four points are shown for $S^j(j)$ with $q=1$ (right figure). This is due to the fact that $S^j(j)$ converges to the exact solution for $j=5$ and the error of $S^5(5)$ is 0. Lines are a guide to the eye.}
\end{figure}
\end{center}

\clearpage

                                                \section{Applications}  
As a first application, let us consider the Taylor expansion of ln$(1+q)$ for small $q$, truncated after $n+1$ terms:

\begin{equation}\label{Taylor_ln}
\text{ln}(1+q) \approx \sum\limits_{k=0}^{n} \frac{(-1)_k}{k!}\frac{q^{k+1}}{k+1}
\end{equation}

which, for $q<1$ can be formally derived from the application of the binomial expansion to the integrand in:

\begin{equation}\label{Taylor_integral}
\text{ln}(1+q) = \int_0^q \frac{1}{1+x}\mathrm{d}x
\end{equation}

followed by the integration. Rather than applying the standard binomial expansion, if one uses the linear transformation in Eq.~\ref{modified_expansion} to expand the integrand, one obtains:

\begin{equation}\label{Linear_transf_ln}
\text{ln}(1+q) \approx \sum\limits_{k=0}^{2j} c_{kj}\frac{(-1)_k}{k!}\frac{q^{k+1}}{k+1}
\end{equation}

The approximations given by Eq.~\ref{Taylor_ln} and Eq.~\ref{Linear_transf_ln} are compared in Figure 4. In both cases the expansion was truncated after 5 terms, i.e. $n=4$ and $j=2$ have been used.
From the figure it can be observed that the Taylor expansion starts diverging in proximity of $q=1$, while the use of Eq.~\ref{Linear_transf_ln} provides a good
approximation of the logarithmic function in an broader interval, well beyond $q=1$.

\begin{center}
\begin{figure}
\includegraphics[scale=0.65]{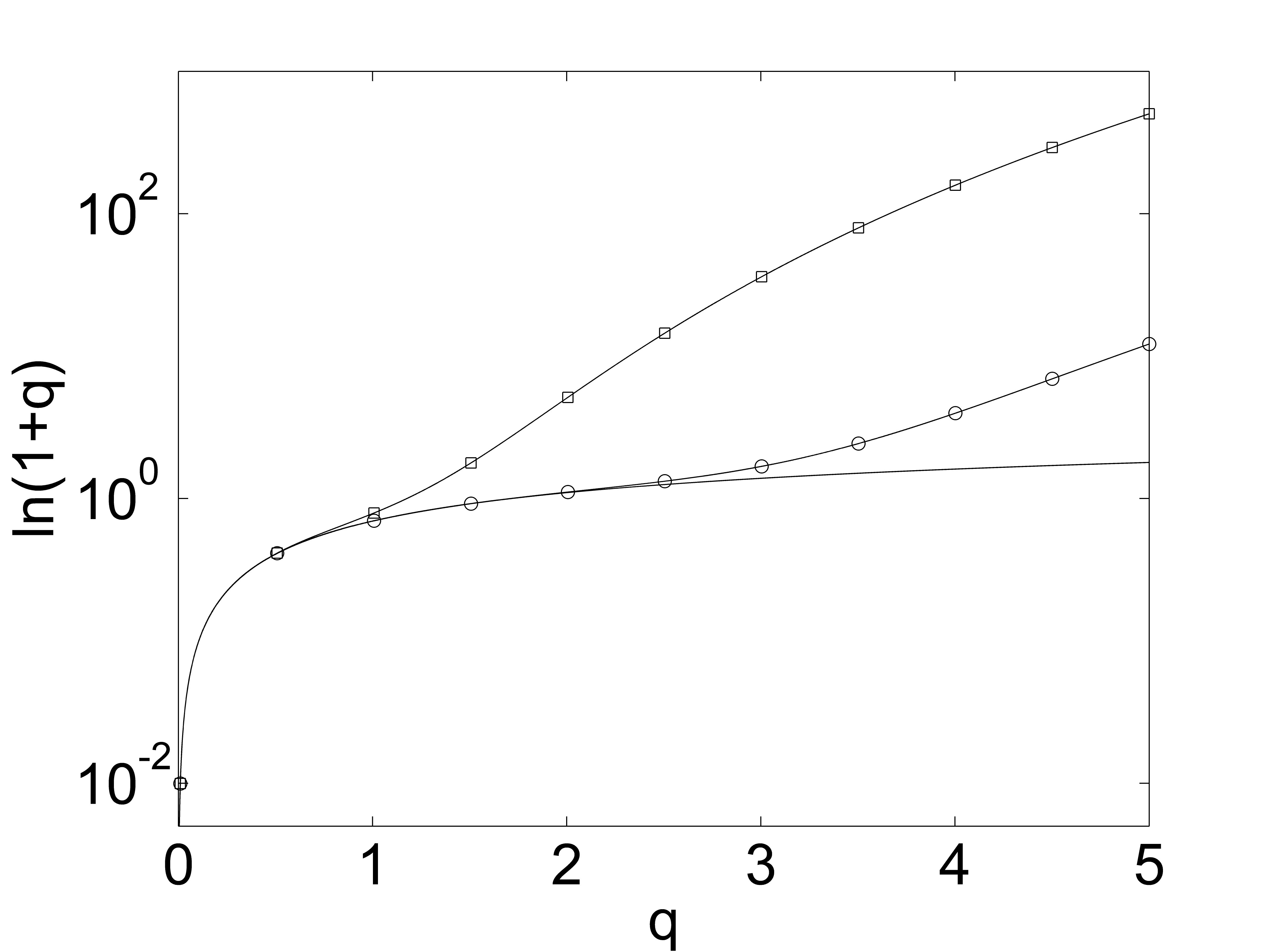}
\caption{\label{Taylor_of_ln_with_5_terms} Plot of ln$(1+q)$. Exact function (line with no symbols), Taylor approximation (squares), approximation using Eq.~\ref{Linear_transf_ln} (circles). Both polynomial approximations are computed with series truncated after 5 terms.}
\end{figure}
\end{center}

As a second application, the computation of the incomplete beta function is considered. For $a,b>0$ and $0\leq x \leq 1$, the incomplete beta function is defined as \cite{Walck}:

\begin{equation}\label{beta_function}
B_x(a,b) = \int_0^x t^{a-1}(1-t)^{b-1}\mathrm{d}t
\end{equation}

Several alternative methods are available for approximating the integral in Eq.~\ref{beta_function} \cite{Walck,Flannery,Stegun}. Binomial expansion of $(1-t)^{b-1}$ followed by integration leads to:
\begin{equation}\label{beta_binomial}
B_x(a,b) \approx x^a\sum\limits_{k=0}^{n} \frac{(-1)^k(b-1)_k}{k!}\frac{x^{k}}{k+a}
\end{equation}

where we have again truncated the series after $n+1$ terms.
Alternatively, one can approximate the incomplete beta function by means of continued fractions:
\begin{equation}\label{beta_cont_fr}
B_x(a,b) \approx \frac{x^a(1-x)^b}{a}\frac{1}{1+}\frac{d_1}{1+}\frac{d_2}{1+}...
\end{equation}

where
\begin{equation}\label{d2mp1_cont_fr}
d_{2m+1}= - \frac{(a+m)(a+b+m)x}{(a+2m)(a+2m+1)}
\end{equation}

and
\begin{equation}\label{d2m_cont_fr}
d_{2m}= \frac{m(b-m)x}{(a+2m-1)(a+2m)}
\end{equation}

An alternative polynomial expansion can be derived by applying the herein proposed transformation. For this, let us consider the alternative definition of $B_x(a,b)$ \cite{Walck}:
\begin{equation}\label{beta_function2}
B_x(a,b) = \int_0^{\frac{x}{1-x}} \frac{u^{a-1}}{(1+u)^{a+b}}\mathrm{d}u
\end{equation}

By applying Eq.~\ref{modified_expansion} to $(1+u)^{-a-b}$, followed by integration, one obtains:
\begin{equation}\label{beta_lin_transf}
B_x(a,b) \approx \sum\limits_{k=0}^{2j} c_{kj}\frac{(-a-b)_k}{k!}\frac{u^{a+k}}{k+a}   \qquad \text{for } u\leq1 \text{, with } u=\frac{x}{1-x}
\end{equation}

and

\begin{equation}\label{beta_lin_transf2}
B_x(a,b) \approx \sum\limits_{k=0}^{2j} c_{kj}\frac{(-a-b)_k}{k!} \left(\frac{1}{k+a}+\frac{1}{k+b}-\frac{1}{(k+b)u^{b+k}}\right)   \qquad \text{for } u>1 \text{, with } u=\frac{x}{1-x}
\end{equation}

A comparison between the three expressions, Eq.~\ref{beta_binomial}, Eq.~\ref{beta_cont_fr} and Eqs.~\ref{beta_lin_transf},~\ref{beta_lin_transf2} is shown in Figure 5 for two non-integer
values of the parameters $a$ and $b$.
All three forms approximate well the incomplete beta function over a wide range of $x$ values, but the accuracy of both Eq.~\ref{beta_binomial} and Eq.~\ref{beta_cont_fr} deteriorates
significantly as $x$ approaches 1, while the beta function calculated by means of the linear transformation introduced in this work remains practically not distinguishable from the correct solution. A careful analysis of the error, shown in Figure 6, evidences that Eq.~\ref{beta_cont_fr} is definitely the most accurate
over a wide range of $x$ values, outperforming the other two approximations for small $x$. On the other hand, from the same figure it is evident that Eqs.~\ref{beta_lin_transf} and ~\ref{beta_lin_transf2}, despite being the less accurate for small $x$, can be applied over the whole range of $x$, which represents a significant advantage over the other two approximations.

Following a similar approach, an accurate polynomial expansion can also be derived for the Student's $t$-distribution, as described in a separate paper \cite{Costa2016}.

\begin{center}
\begin{figure}
\includegraphics[scale=0.65]{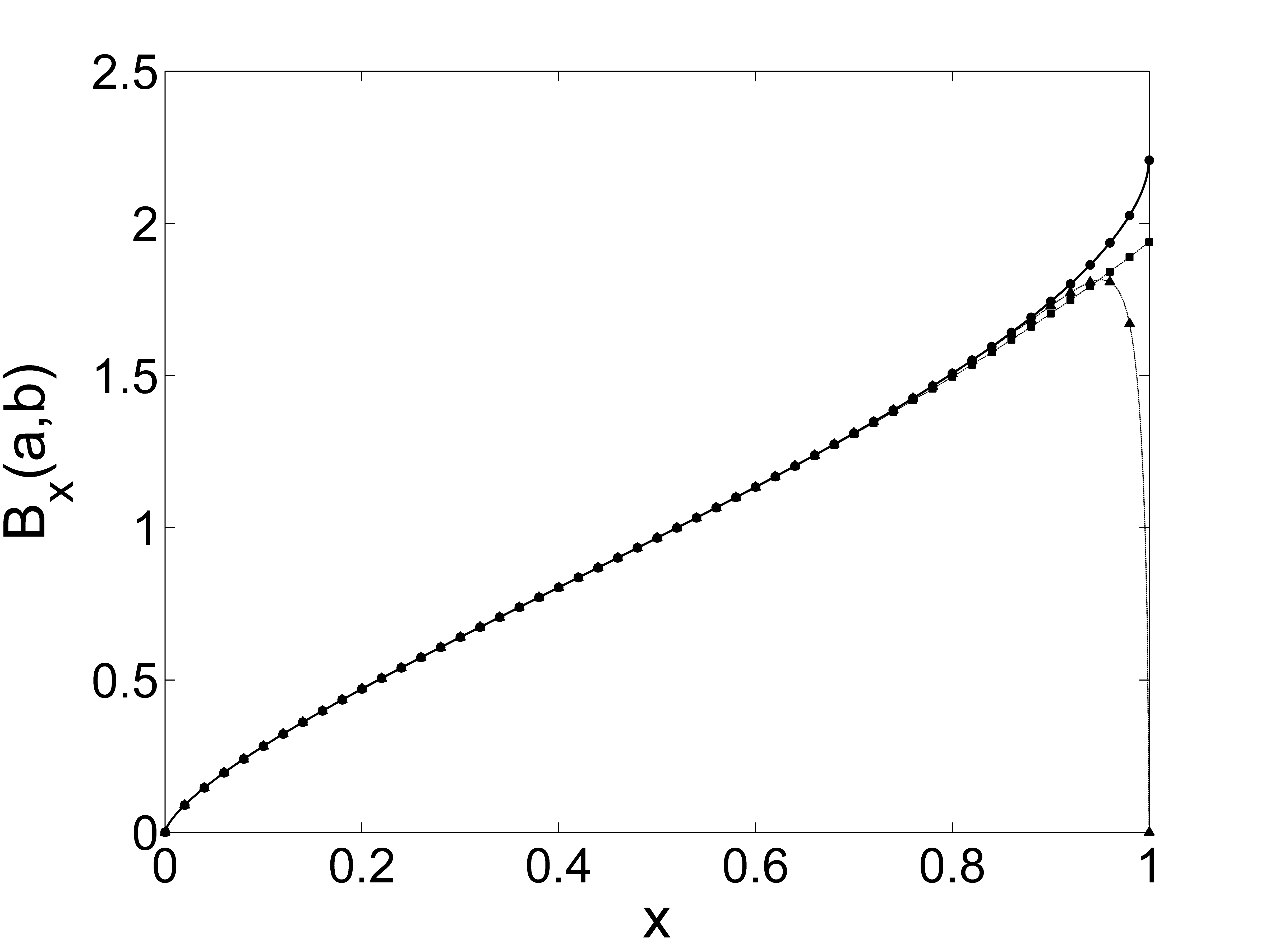}
\caption{\label{Beta_function} Plot of $B_x(a,b)$ for $a = 1/\sqrt[]{2}$ and $b = 1/\sqrt[]{3}$: exact (bold continuous line), calculated with Eq.~\ref{beta_binomial} (squares), calculated with
 Eq.~\ref{beta_cont_fr} (triangles), calculated with Eqs.~\ref{beta_lin_transf},~\ref{beta_lin_transf2} (circles). The three approximated solutions were all computed using 7 terms.}
\end{figure}
\end{center}

\begin{center}
\begin{figure}
\includegraphics[scale=0.65]{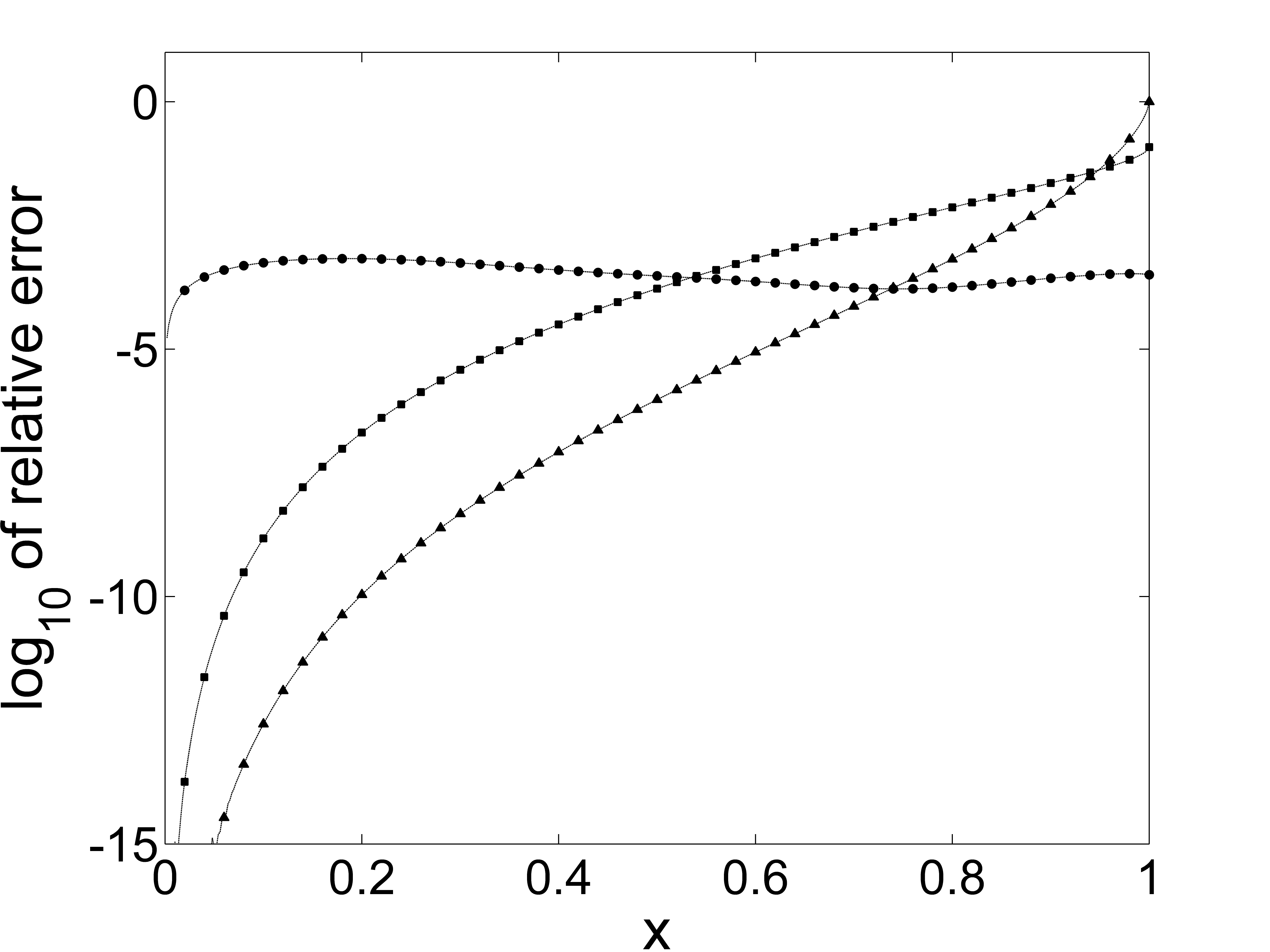}
\caption{\label{Beta_function_err} Plot of the relative error in calculating $B_x(a,b)$ for $a = 1/\sqrt[]{2}$ and $b = 1/\sqrt[]{3}$ using 7 terms: Eq.~\ref{beta_binomial} (squares),
 Eq.~\ref{beta_cont_fr} (triangles), Eqs.~\ref{beta_lin_transf},~\ref{beta_lin_transf2} (circles).}
\end{figure}
\end{center}

                                                \section{Conclusions}  

To summarize, a simple linear transformation has been proposed to accelerate the convergence of the binomial series for the case of
negative exponents and binomial terms of the same sign. Several numerical examples have been reported, indicating that the acceleration is significant when the terms of the binomial are of comparable magnitude. Quite remarkably, the transformed series converges also when the terms are equal, i.e. when the Binomial Theorem is not valid.
This allows to derive accurate polynomial expansions of ln$(1+x)$ and of the incomplete beta function whit a broader range of applicability than those
obtained through the Binomial Theorem.
It is worth mentioning that the examples reported in this work considered exponent values ranging from -0.5 to -10. Further analyses
are therefore required to verify the scaling behavior of the truncation error when even smaller exponents are considered.

\end{document}